\documentclass[smallextended]{svjour3}       
\usepackage{graphicx}
\usepackage{epstopdf}
\usepackage{amsmath}
\usepackage{amssymb}
\usepackage{fancyhdr}
\usepackage{color,url}
\usepackage{algorithm,algorithmic}
\usepackage{subfigure}
\usepackage{multirow}
\usepackage{rotating}
\usepackage{tabularx}
\usepackage{booktabs}
\usepackage{lscape}
\allowdisplaybreaks[4]

\begin{document}

\title{A note on $R$-linear convergence of nonmonotone gradient methods
}

\titlerunning{A note on $R$-linear convergence of nonmonotone gradient methods}        

\authorrunning{Yakui Huang, Yu-Hong Dai, Xin-Wei Liu} 

\author{Xinrui Li         \and
	Yakui Huang 
}


\institute{
	Xinrui Li 
	\and
	Yakui Huang \at
	Institute of Mathematics, Hebei University of Technology, Tianjin 300401, China\\
	\email{hyk@hebut.edu.cn}
}

\date{Received: date / Accepted: date}

\maketitle

\begin{abstract}
Nonmonotone gradient methods generally perform better than their monotone counterparts especially on unconstrained quadratic optimization. However, the known convergence rate of the monotone method is often much better than its nonmonotone variant. With the aim of shrinking the gap between theory and practice of nonmonotone gradient methods, we introduce a property for convergence analysis of a large collection of gradient methods. We prove that any gradient method using stepsizes satisfying the property will converge $R$-linearly at a rate of $1-\lambda_1/M_1$, where $\lambda_1$ is the smallest eigenvalue of Hessian matrix and $M_1$ is the upper bound of the inverse stepsize. Our results indicate that the existing convergence rates of many nonmonotone methods can be improved to $1-1/\kappa$ with $\kappa$ being the associated condition number.	
\end{abstract}

\section{Introduction}\label{intro}
We concentrate on analyzing the convergence rate of gradient methods for unconstrained quadratic optimization:
\begin{equation}\label{1a}
	\min_{x\in \mathbb{R}^{n}} f(x)=\frac{1}{2}x^{\top}Ax-b^{\top}x,
\end{equation}
where $A\in \mathbb{R}^{n\times n}$ is a real symmetric positive definite  matrix and $b\in \mathbb{R}^{n}$. The gradient method solves \eqref{1a} by
\begin{equation}\label{aa}
	x_{k+1}=x_{k}-\alpha_{k}g_{k},
\end{equation}
where $g_{k} = \nabla f(x_{k})$ and  $\alpha_{k}>0$ is the stepsize. By considering the generated function values, we can roughly divided gradient methods into two categories: monotone and nonmonotone.

The classic steepest descent (SD) method developed by Cauchy \cite{cauchy1847methode} is monotone, which determines the stepsize by exact line search:
\begin{equation}\label{sd}
	\alpha_k^{SD}=\arg\min_{\alpha>0}~f(x_k-\alpha g_k)=\frac{g_k^Tg_k}{g_k^TAg_k}.
\end{equation}
It is known that the SD method converges $Q$-linearly  at a rate of $\frac{\kappa-1}{\kappa+1}$ \cite{akaike1959successive,forsythe1968asymptotic}, where $\kappa=\lambda_{n}/\lambda_{1}$ is the condition number of $A$ with $\lambda_{1}$ and $\lambda_{n}$ being the smallest and largest eigenvalues, respectively. In a similar spirit to the SD method, the minimal gradient (MG) method  \cite{dai2003altermin} calculates the stepsize by minimizing the gradient norm, i.e.
\begin{equation}\label{mg}
	\alpha_k^{MG}=\arg\min_{\alpha>0}~\|g(x_{k}-\alpha g_{k})\|=\frac{g_k^TAg_k}{g_k^TA^2g_k}.
\end{equation}
This method is also monotone and has the same convergence speed as the SD method \cite{huang2019asymptotic}.
Another monotone method is due to Dai and Yang \cite{dai2006new2} whose stepsize is given by
\begin{equation}\label{saopt}
	\alpha_k^{AOPT}=\frac{\|g_k\|}{\|Ag_k\|}=\sqrt{\alpha_k^{SD}\alpha_k^{MG}}.
\end{equation}
It is proved that $\alpha_k^{AOPT}$ asymptotically converges to $\frac{2}{\lambda_1+\lambda_n}$, which is in some sense an optimal stepsize since it minimizes $\|I-\alpha A\|$ over $\alpha$, see \cite{dai2006new2,elman1994inexact}. The Dai--Yang method \eqref{saopt} also converges $Q$-linearly at the same rate as the SD method \cite{dai2006new2}. See \cite{dai2003altermin,dai2005analysis,yuan2006new,yuan2008step} and references therein for other monotone gradient methods.

The Barzilai--Borwein (BB) methods \cite{barzilai1988two} can be viewed as nonmonotone counterparts of the SD and MG methods since
\begin{equation}\label{bb1}
	\alpha_{k}^{BB1}=\frac{s_{k-1}^{\top} s_{k-1}}{s_{k-1}^{\top} y_{k-1}}=\frac{g_{k-1}^{\top} g_{k-1}}{g_{k-1}^{\top} A g_{k-1}}=\alpha_{k-1}^{SD}
\end{equation}
and
\begin{equation}\label{bb2}
	\alpha_{k}^{BB2}=\frac{s_{k-1}^{\top} y_{k-1}}{y_{k-1}^{\top} y_{k-1}}=\frac{g_{k-1}^{\top} A g_{k-1}}{g_{k-1}^{\top} A^{2} g_{k-1}}=\alpha_{k-1}^{MG},
\end{equation}
where $s_{k-1}=x_{k}-x_{k-1}$ and $y_{k-1}=g_{k}-g_{k-1}$. Clearly, $\alpha_{k}^{BB1} \geqslant  \alpha_{k}^{BB2}$ for the quadratic function \eqref{1a}. An interesting property of the above two BB stepsizes is that they satisfy the quasi-Newton secant equation in the sense of least squares:
\begin{equation*}\label{bbquasi}
	\alpha_{k}^{BB1}=\arg \min _{\alpha>0}\left\|\alpha^{-1} s_{k-1}-y_{k-1}\right\|,\quad
	\alpha_{k}^{BB2}=\arg \min _{\alpha>0}\left\|s_{k-1}-\alpha y_{k-1}\right\|.
\end{equation*}
For the BB2 method, a surprising $R$-superlinear convergence result on the two dimensional strictly convex quadratic function  is presented in \cite{barzilai1988two}, which is improved by Dai \cite{dai2013new} for the BB1 method. Dai and Fletcher \cite{dai2005asymptotic} showed similar $R$-superlinear convergence for the BB1 method in the three dimensional case. As for the $n$ dimensional case, Raydan \cite{raydan1993barzilai} established global convergence of the BB1 method while Dai and Liao \cite{dai2002r} proved that it converges $R$-linearly and the rate is roughly $(1-1/\kappa)^{1/n}$. Recently, Li and Sun \cite{li2021on} improved the above rate to $1-1/\kappa$. A nonmonotone version of the Dai--Yang method using $\alpha_{k-1}^{AOPT}$ can be found in \cite{dai2015positive} which is also $R$-superlinear convergent on two dimensional strictly convex quadratics.

Interestingly, the above mentioned nonmonotone methods are generally much faster than the corresponding monotone ones especially for quadratics. Moreover, those nonmonotone methods have a great
advantage of easy extension to general unconstrained problems \cite{BDH2019stabilized,dai2019family,di2018steplength,fletcher2005barzilai,huang2020equipping,raydan1997barzilai}, to constrained optimization \cite{birgin2000nonmonotone,crisci2020spectral,dai2005asymptotic,dai2005projected,huang2020equipping,huang2016smoothing} and to various applications \cite{dai2015positive,huang2020gradient,huang2015quadratic,jiang2013feasible,tan2016barzilai,wang2007projected}. However, theoretical analysis of a nonmonotone method is generally much more difficult due to the nonmonotonity. In recent years, to achieve good performance, long and short stepsizes are often employed in alternate, cyclic, adaptive, or hybrid way by nonmonotone gradient methods  \cite{dai2005projected,de2014efficient,frassoldati2008new,huang2020equipping,sun2020new,zhou2006gradient} which makes the analysis even  more difficult. Consequently, the known convergence rate of a nonmonotone method is often much worse than its monotone counterpart.


An important tool for analyzing the convergence of nonmonotone methods is the following Property A suggested by Dai \cite{dai2003alternate}.

\textbf{Property A:} Assume that $A=\operatorname{diag}\left(\lambda_{1}, \cdots, \lambda_{n}\right)$ and $1 =\lambda_{1}<\lambda_{2}<\cdots<\lambda_{n}$. Let $g_{k}^{(i)}$ be the $i$-th component of $g_{k}$ and
$$
P(k, l)=\sum_{i=1}^{l}\left( g_{k}^{(i)}\right)^{2}.
$$
Suppose that there exist an integer $m$ and positive constants  $M_{1} \geqslant  \lambda_{1}$ and $M_{2}$ such that\\
(i) $\lambda_{1} \leqslant  \alpha_{k}^{-1} \leqslant  M_{1}$\\
(ii) for any integer $l \in[1, n-1]$ and $\varepsilon>0$, if $P(k-j, l) \leqslant  \varepsilon$ and $\left(g_{k-j}^{(l+1)}\right)^{2} \geqslant  M_{2} \varepsilon$ hold for $j \in[0, \min \{k, m\}-1]$, then $\alpha_{k}^{-1} \geqslant  \frac{2}{3} \lambda_{l+1}$. 

Dai \cite{dai2003alternate} proved that if stepsizes of the gradient method have Property A, then either $g_{k} = 0$ for some finite $k$ or the sequence  $\left\{\left\|g_{k}\right\|\right\}$ converges to zero $R$-linearly and  the rate is roughly $\big(\max\{\frac{1}{2},1-\frac{\lambda_{1}}{M_1}\}\big)^{1/n}$. Such a rate has recently been improved to  $\max\{\frac{1}{2},1-\frac{\lambda_{1}}{M_1}\}$ \cite{huang2021On}, and the rate can be further improved to $1-1/\kappa$ for the gradient method with:
\begin{equation}\label{gmr}
	\alpha_{k}^{GMR}=\frac{g_{v(k)}^{\top}A^{\rho(k)}g_{v(k)}}{g_{v(k)}^{\top}A^{\rho(k)+1}g_{v(k)}},
\end{equation}
where $v(k)\in\{k,k-1,\dots,\max\{k-m+1,1\}\}$ and $\rho(k)\geqslant 0$ is some real number, see \cite{friedlander1998gradient}. Nevertheless, it is not very straightfoward to check Property A since it requires that for any integer $l$, the inverse stepsize is large enough when the first $l$ eigencomponents of $g_{k}$ are much smaller than the $(l + 1)$-th eigencomponent. Moreover, as we will see in the next section some existing stepsize does not satisfy Property A.

In this paper, we introduce a new property, named Property B, for stepsizes of the gradient method which is straightforward and easy to check. We prove that any gradient method with stepsizes satisfying  Property B will converge $R$-linearly at a rate of $1-\lambda_{1}/M_1$ with $M_1$ being the upper bound of the inverse stepsize as above. Our results indicate that gradient methods using retard, alternate, cyclic, adaptive or hybrid stepsize schemes such as the  SDC   method \cite{de2014efficient}, the periodic gradient method in \cite{huang2019asymptotic} and the methods in \cite{huang2020acceleration,sun2020new,zou2021fast} converge $R$-linearly at the above rate. Consequently, when $M_1\leqslant \lambda_{n}$ the convergence rate can be refined to $1-1/\kappa$ which implies that the BB1, BB2, alternate BB \cite{dai2005projected}, adaptive BB \cite{zhou2006gradient}, cyclic BB \cite{dai2006cyclic}, cyclic SD \cite{dai2005asymptotic}, ABBmin \cite{frassoldati2008new}, BB family \cite{dai2019family}, and the methods in \cite{dai2015positive,huang2020equipping,huang2019asymptotic,huang2020gradient,huang2020acceleration} converge at least $R$-linearly with a rate of $1-1/\kappa$. It is worth to mention that the methods in \cite{huang2019asymptotic,huang2020acceleration} satisfy Property B but not necessarily satisfy Property A. In addition, our results generalize the one in \cite{li2021on} for the BB1 method, and improve those in \cite{huang2021On} for the method \eqref{gmr} in the sense that the associated quantities in our convergence rate are smaller than or equal to those in \cite{huang2021On}. We also discuss convergence rates of iterations and function values for gradient methods satisfying the proposed Property B.

The paper is organized as follows. In Section \ref{sec:2},  we present Property B and our $R$-linear convergence results. In Section \ref{sec:3} we give some concluding remarks and discuss generalizations of our results.

\section{Main results}\label{sec:2}
Since the gradient method \eqref{aa} is invariant under translations and rotations when applying to the quadratic problem \eqref{1a}, for the convenience of analysis, we assume without loss of generality that the matrix $A$ is diagonal  with distinct eigenvalues, i.e.
\begin{equation}\label{diag}
	A=\mathrm{diag}\{\lambda_{1}, \lambda_{2}, \cdots, \lambda_{n}\}, \quad 0<\lambda_{1}<\lambda_{2}<\cdots<\lambda_{n}.
\end{equation}

We shall analyze $R$-linear convergence rate of gradient methods whose stepsizes satisfying the following Property B.

\noindent\textbf{Property~B:} We
say that the stepsize $\alpha_{k}$ has Property B if there exist  an integer $m$ and a positive constant  $M_{1} \geqslant  \lambda_{1}$ such that\\
(i) $\lambda_{1} \leqslant  \alpha_{k}^{-1} \leqslant  M_{1}$;\\
(ii) for some integer $v(k) \in\{k,k-1,\cdots,\max\{k-m+1,0\}\}$,
\begin{equation}\label{alphack} 
	\alpha_{k} \leqslant  \frac{g_{v(k)}^{\top} \psi(A) g_{v(k)}}{g_{v(k)}^{\top} A\psi(A) g_{v(k)}},
\end{equation}
where $\psi$ is a real analytic function  on $[\lambda_{1},\lambda_{n}]$ and  can be expressed by Laurent series 
\begin{equation*}\label{psi}
	\psi(z)=\sum_{k=-\infty}^{\infty} d_{k} z^{k},\quad d_{k} \in \mathbb{R},
\end{equation*}
such that $0<\sum_{k=-\infty}^{\infty} d_{k} z^{k}<+\infty$ for all $z \in\left[\lambda_{1}, \lambda_{n}\right]$. 

The only difference between Property B here and Property A is condition (ii). For the case $\psi(A) = I$, we can prove that $\alpha_{k}$ satisfies \eqref{alphack} also meets the condition (ii) of Property A with $M_{2} = 2$. In fact, for any integer $ l \in [1,n-1] $, suppose that $P(k-j, l) \leqslant  \varepsilon$ and $\left(g_{k-j}^{(l+1)}\right)^{2} \geqslant  2 \varepsilon$ hold for $\varepsilon>0$ and $j \in[0, \min \{k, m\}-1]$, we obtain that
\begin{equation*}
	\begin{aligned}
		\alpha_{k}^{-1} & \geqslant  \frac{g_{v(k)}^{\top} A g_{v(k)}}{g_{v(k)}^{\top}  g_{v(k)}}
		=\frac{\sum_{i=1}^{n} \lambda_{i}\left(g_{v(k)}^{(i)}\right)^{2}}{\sum_{i=1}^{n}\left(g_{v(k)}^{(i)}\right)^{2}} \geqslant  \frac{\lambda_{l+1} \sum_{i=l+1}^{n}\left(g_{v(k)}^{(i)}\right)^{2}}{P(v(k), l)+\sum_{i=l+1}^{n}\left(g_{v(k)}^{(i)}\right)^{2}} \\
		& \geqslant  \frac{\lambda_{l+1}}{\varepsilon / 2 \varepsilon+1} \geqslant  \frac{2}{3} \lambda_{l+1}.
	\end{aligned}
\end{equation*}
However, if $\psi(A) \neq I$, the following example shows that $\alpha_{k}^{-1}<\frac{2}{3} \lambda_{l+1}$ for a stepsize satisfies \eqref{alphack}  even if $P(k-j, l) \leqslant  \varepsilon$ and $\left(g_{k-j}^{(l+1)}\right)^{2} \geqslant  M_2 \varepsilon$.


\textbf{Example 1}. Consider the $3$-dimensional quadratic with 
\begin{equation*}
	A=\mathrm{diag}\{1, 8, 16\},\quad b=0.
\end{equation*}
We apply the gradient method in \cite{huang2020acceleration} which uses the stepsize 
\begin{equation*}
	\alpha_{k}
	=\frac{g_{k-1}^{\top} \psi(A) g_{k-1}}{g_{k-1}^{\top} A\psi(A) g_{k-1}},~k\geqslant 1.
\end{equation*}
Here we choose $\psi(z)=\left(\frac{1+2z}{z^2}\right)^2$ and $\alpha_0 = \alpha_0^{SD}$. For any given  $\varepsilon>0$, set   $x_0=\left(\sqrt{\varepsilon},\frac{\sqrt{40\varepsilon}}{8},\frac{\sqrt{40\varepsilon}}{16}\right)^T$. We get
\begin{equation*}
	g_0=Ax_0=\left(\sqrt{\varepsilon},\sqrt{40\varepsilon},\sqrt{40\varepsilon}\right)^T.
\end{equation*}	
It follows that $\alpha_0=\frac{g_{0}^{\top} g_{0}}{g_{0}^{\top} Ag_{0}}\approx0.084~3$ and
\begin{equation*}
	g_1=(I-\alpha_0A)g_0
	\approx\left(0.915~7\sqrt{\varepsilon},2.059~9\sqrt{\varepsilon},-2.204~7\sqrt{\varepsilon}\right)^T,
\end{equation*}
which gives $\alpha_{1}=\frac{g_{0}^{\top} \psi^{2}(A) g_{0}}{g_{0}^{\top} A\psi^{2}(A) g_{0}}\approx0.295~8$
and
\begin{equation*}
	g_2=(I-\alpha_1A)g_1
	\approx\left(0.644~8\sqrt{\varepsilon},-2.814~8\sqrt{\varepsilon},8.230~0\sqrt{\varepsilon}\right)^T.
\end{equation*}
Then we have $\alpha_{2}\approx0.705~6$ and
\begin{equation*}
	g_3=(I-\alpha_2A)g_2
	\approx\left(0.189~8\sqrt{\varepsilon},13.0744\sqrt{\varepsilon},-84.684~6\sqrt{\varepsilon}\right)^T.
\end{equation*}		
Clearly, for $j \in[0, \min \{k, m\}-1]$ with $m=2$, $P(k-j,1) \leqslant  \varepsilon$ holds for  $k=2,3$. However, for any $M_2$ such that $\left(g_{k-j}^{(2)}\right)^{2} \geqslant  M_{2} \varepsilon$, we have $\alpha_{2}^{-1}\approx 1.417~2<16/3$ and $\alpha_{3}^{-1}\approx 4.832~0<16/3$.

Now we present our first convergence result for gradient methods satisfying Property B.	See \cite{li2021on} for analysis of the BB1 method and \cite{huang2021On} for methods satisfying Property A.
\begin{theorem}\label{converth}
	Suppose that the sequence $\left\{\left\|g_{k}\right\|\right\}$ is generated by applying the gradient method \eqref{aa} to the quadratic problem \eqref{1a}. If the stepsize $\alpha_{k}$ satisfies Property B, then either $g_{k} = 0$ for some finite $k$ or the sequence  $\left\{\left\|g_{k}\right\|\right\}$ converges to zero $R$-linearly in the sense that
	\begin{equation}\label{Citheta}
		\lvert g_{k}^{(i)}\rvert \leqslant  C_{i}\theta^{k}, \qquad i=1,2, \cdots, n,
	\end{equation}
	where  ${\theta}=1-\lambda_{1}/M_1$ and
	\begin{equation*}
		\left\{\begin{array}{l}
			C_{1}=\lvert g_{0}^{(1)}\rvert; \\
			C_{i}=\displaystyle\max\left\{\lvert g_{0}^{(i)}\rvert,\frac{\lvert g_{1}^{(i)}\rvert}{\theta}, \cdots,\frac{\lvert g _{m-1}^{(i)}\rvert}{\theta^{m-1}},\frac{\max \{\sigma_{i}, \sigma_{i}^{m}\}}{\theta^{m}\psi(\lambda_{i})}\sqrt{\sum_{j=1}^{i-1}\psi^{2}(\lambda_{j}) C_{j}^{2}}\right\},\\
			\quad i=2,3, \cdots, n,
		\end{array}\right.
	\end{equation*}
	with $\sigma_{i}=\max \left\{\frac{\lambda_{i}}{\lambda_{1}}-1,1-\frac{\lambda_{i}}{M_{1}}\right\}$.
\end{theorem}
\begin{proof}
	Clearly, we only need to show that \eqref{Citheta} holds for the case $g_{k}^{(i)}\neq0$.  In what follows, we show \eqref{Citheta} by induction on $i$. Since $g_{k+1}^{(i)}=\left(I-\alpha_{k}\lambda_{i}\right)g_{k}^{(i)}$ and $\lambda_{1} \leqslant  \alpha_{k}^{-1} \leqslant  M_{1}$,  we have
	\begin{equation}\label{sigmai}
		\lvert g_{k+1}^{(i)}\rvert \leqslant  {\sigma}_{i}\lvert g_{k}^{(i)}\rvert,\qquad i=1,2, \cdots, n,
	\end{equation}
	which gives
	\begin{equation*}\label{eqCitheta1}
		\lvert g_{k}^{(1)}\rvert \leqslant  \theta\lvert g_{k-1}^{(1)}\rvert \leqslant  C_{1}\theta^{k}.
	\end{equation*}
	That is, \eqref{Citheta} holds for $i=1$.

	Assume that \eqref{Citheta} holds for all $1 \leqslant  i \leqslant  L-1$ with $L\in\{2, \cdots, n\}$. We have to prove that \eqref{Citheta} holds for $i=L$. By the definition of  $C_{L}$, we know that \eqref{Citheta} holds for $k=0,1, \cdots, m-1$. Assume by contradiction that \eqref{Citheta} does not hold for $k\geqslant  m$ when $i=L$. Let $\hat{k}~(\geqslant  m)$ be the minimal index such that $\lvert g_{\hat{k}}^{(L)}\rvert>C_{L} \theta^{\hat{k}}$. Then, we must have $\alpha_{{}\hat{k}-1} \lambda_{L}>1$; otherwise,
	\begin{equation*}\label{ineqcontr1}
		\lvert g_{{}\hat{k}}^{(L)}\rvert =
		(1-\alpha_{{}\hat{k}-1} \lambda_{L})\lvert g_{{}\hat{k}-1}^{(L)}\rvert\leqslant  {\theta}\lvert g_{{}\hat{k}-1}^{(L)}\rvert\leqslant  C_{L} \theta^{{}\hat{k}},
	\end{equation*}
	which is impossible.  It follows from \eqref{sigmai}, $\theta<1$ and the inductive assumption that
	\begin{align*}
		\psi(\lambda_{L})\lvert g_{\hat{k}-j}^{(L)}\rvert 
		&\geqslant  \frac{\psi(\lambda_{L})\lvert g_{\hat{k}}^{(L)}\rvert}{{\sigma}_{L}^{j}}>\frac{\psi(\lambda_{L})C_{L} \theta^{\hat{k}}}{{\sigma}_{L}^{j}} 
		\geqslant  \frac{\psi(\lambda_{L})C_{L} \theta^{\hat{k}}}{\max \{\sigma_{L}, \sigma_{L}^{m}\}} 
		\\
		& \geqslant  \theta^{\hat{k}-m}\sqrt{\sum_{i=1}^{L-1}\psi^{2}(\lambda_{i}) C_{i}^{2}} \nonumber
		\geqslant  \theta^{\hat{k}-m}\sqrt{ \sum_{i=1}^{L-1}\left(\psi(\lambda_{i})g_{\hat{k}-j}^{(i)}\right)^{2} \theta^{2(j-\hat{k})}}  \nonumber\\
		&=\theta^{j-m} \sqrt{G(\hat{k}-j, L-1)} 
		\geqslant  \sqrt{G(\hat{k}-j, L-1)},\quad j \in[1, m],
	\end{align*}
	which gives
	\begin{equation}\label{ineqGeps2}
		\left(\psi(\lambda_{L})g_{v(\hat{k}-1)}^{(L)}\right)^2> G(v(\hat{k}-1), L-1),
	\end{equation}
	where $G(k, l)=\sum_{i=1}^{l}\left(\psi(\lambda_{i})g_{k}^{(i)}\right)^{2}$. Since $\alpha_{\hat{k}-1} \lambda_{L}>1$, by \eqref{alphack} and \eqref{ineqGeps2} we have that
	\begin{align}\label{ineqlin1}
		\left\lvert1-\alpha_{\hat{k}-1} \lambda_{L}\right\rvert&=\alpha_{\hat{k}-1} \lambda_{L}-1
		\leqslant  \frac{\left(\psi(A) g_{v(\hat{k}-1)}\right)^{\top}\left(\psi(A) g_{v(\hat{k}-1)}\right)}{\left(\psi(A) g_{v(\hat{k}-1)}\right)^{\top} A\left(\psi(A) g_{v(\hat{k}-1)}\right)}\lambda_{L}-1 \nonumber\\
		&=\frac{\sum_{i=1}^{n}\left(\lambda_{L}-\lambda_{i}\right) \left(\psi(\lambda_{i})g_{v(\hat{k}-1)}^{(i)}\right)^{2}}{\sum_{i=1}^{n} \lambda_{i}\left(\psi(\lambda_{i})g_{v(\hat{k}-1)}^{(i)}\right)^{2}}\nonumber\\ 
		&\leqslant  \frac{\left(\lambda_{L}-\lambda_{1}\right)G(v(\hat{k}-1), L-1)}{\lambda_{1}G(v(\hat{k}-1), L-1) +\lambda_{L}\left(\psi(\lambda_{L})g_{v(\hat{k}-1)}^{(L)}\right)^{2} } \nonumber\\	
		&\leqslant \frac{\lambda_{L}-\lambda_{1}}{\lambda_{1}+\lambda_{L}\frac{ \left(\psi(\lambda_{L})g_{v(\hat{k}-1)}^{(L)}\right)^{2}}{ G(v(\hat{k}-1), L-1)}} \leqslant  \frac{\lambda_{L}-\lambda_{1}}{\lambda_{L}+\lambda_{1}}.
	\end{align}
	Notice that	
	\begin{align*}\label{ineqlin2}
		\alpha_{\hat{k}-1}^{-1}
		&\geqslant \frac{\sum_{i=1}^{n}\lambda_{i} \left(\psi(\lambda_{i})g_{v(\hat{k}-1)}^{(i)}\right)^{2}}{\sum_{i=1}^{n} \left(\psi(\lambda_{i})g_{v(\hat{k}-1)}^{(i)}\right)^{2}}\nonumber\\ 
		&\geqslant  \frac{\lambda_{1}G(v(\hat{k}-1), L-1)+\lambda_{L}\sum_{i=L}^{n} \left(\psi(\lambda_{i})g_{v(\hat{k}-1)}^{(i)}\right)^{2}}{G(v(\hat{k}-1), L-1)+\sum_{i=L}^{n} \left(\psi(\lambda_{i})g_{v(\hat{k}-1)}^{(i)}\right)^{2}} \nonumber\\
		&=\frac{\lambda_{1}+\lambda_{L}\delta}{1+\delta}>\frac{\lambda_{1}+\lambda_{L}}{2},		
	\end{align*}
	where $\delta=\frac{\sum_{i=L}^{n} \left(\psi(\lambda_{i})g_{v(\hat{k}-1)}^{(i)}\right)^{2}}{G(v(\hat{k}-1), L-1)}\geqslant \frac{\left(\psi(\lambda_{L})g_{v(\hat{k}-1)}^{(L)}\right)^2}{G(v(\hat{k}-1), L-1)}>1$. It follows from (i) of Property B that $M_1>\frac{\lambda_{1}+\lambda_{L}}{2}$. Thus,
	\begin{equation*}
		\frac{\lambda_{L}-\lambda_{1}}{\lambda_{L}+\lambda_{1}}=1-\frac{2\lambda_{1}}{\lambda_{L}+\lambda_{1}}<1-\frac{\lambda_{1}}{M_{1}}=\theta,
	\end{equation*}
	which together with	\eqref{ineqlin1} gives $\lvert g_{\hat{k}}^{(L)}\rvert\leqslant  {\theta}\lvert g_{\hat{k}-1}^{(L)}\rvert\leqslant  C_{L} \theta^{\hat{k}}$. This contradicts our assumption. Thus \eqref{Citheta} holds for all $i$. We complete the proof.
\end{proof}

Theorem \ref{converth} allows us to use various choices of stepsizes for the gradient method. For example, one can use a stepsize such that $v(k)=k-r$ for some integer $r\in[1,m-1]$ and all $k\geqslant  r$ which includes BB1, BB2, and the method in \cite{dai2015positive} as special instances. One may also use alternate, cyclic, adaptive or hybrid stepsize schemes such as the cyclic BB \cite{dai2006cyclic}, ABBmin \cite{frassoldati2008new} and
SDC  \cite{de2014efficient} methods, the periodic gradient method in \cite{huang2019asymptotic} and the methods in \cite{huang2020gradient,sun2020new,zou2021fast}.  Convergence results for those cases are specified in the next corollary.
\begin{corollary}\label{converfix1}
	Suppose that the sequence $\left\{\|g_{k}\|\right\}$ is generated by applying the gradient method \eqref{aa} to the quadratic problem \eqref{1a}.
	\begin{enumerate}
		\item[(i)] If Property B is satisfied with $v(k)=k-r$ for some integer $r\in[1,m-1]$ and all $k\geqslant  r$, then either $g_{k} = 0$ for some finite $k$ or the sequence  $\left\{\|g_{k}\|\right\}$ converges to zero $R$-linearly in the sense of \eqref{Citheta} with
		\begin{equation*}
			\left\{\begin{array}{l}
				C_{1}=\lvert g_{0}^{(1)}\rvert; \\
				C_{i}=\displaystyle
				\max\left\{
				\lvert g_{0}^{(i)}\rvert,\frac{\lvert g_{1}^{(i)}\rvert}{\theta}, \cdots,\frac{\lvert g_{r}^{(i)}\rvert}{\theta^{r}},\frac{\sigma_{i}^{r+1}}{\theta^{r+1}\psi(\lambda_{i})}\sqrt{\sum_{j=1}^{i-1}\psi^{2}(\lambda_{j}) C_{j}^{2}}
				\right\},\\
				\quad i=2,3, \cdots, n.
			\end{array}\right.
		\end{equation*}

		\item[(ii)] If Property B is satisfied by replacing $v(k)$ in \eqref{alphack} with $v_i(k)\in\{k,k-1,\cdots,\max\{k-m_i+1,0\}\}$ for all $i=1,\cdots,s$ where $m_i$ and $s$ are positive integers, then either $g_{k} = 0$ for some finite $k$ or the sequence  $\left\{\|g_{k}\|\right\}$ converges to zero $R$-linearly in the sense of \eqref{Citheta} with $m=\max\{m_1,\cdots,m_s\}$.
	\end{enumerate}	
\end{corollary}

Many efficient gradient methods using stepsizes satisfying $\lambda_{1}\leqslant \alpha_{k}^{-1}\leqslant  \lambda_{n}$, for example, the BB1, BB2, alternate BB \cite{dai2005projected}, adaptive BB \cite{zhou2006gradient}, cyclic BB \cite{dai2006cyclic}, and cyclic SD \cite{dai2005asymptotic} methods, BB family \cite{dai2019family}, and the method in \cite{dai2015positive}. For such methods, next theorem refines the results in Theorem \ref{converth} and Corollary \ref{converfix1}.	
\begin{theorem}\label{converth2}
	Suppose that the sequence $\left\{\|g_{k}\|\right\}$ is generated by applying the gradient method \eqref{aa} with $\lambda_{1}\leqslant \alpha_{k}^{-1}\leqslant  \lambda_{n}$ to the quadratic problem \eqref{1a}.
	\begin{enumerate}
		\item[(i)] If Property B is satisfied with $v(k)=k-r$ for some integer $r\in[1,m-1]$ and all $k\geqslant  r$, then either $g_{k} = 0$ for some finite $k$ or the sequence  $\left\{\|g_{k}\|\right\}$ converges to zero $R$-linearly in the sense that
		\begin{equation}\label{Citheta3}
			\lvert g_{k}^{(i)}\rvert \leqslant  \widetilde{C}_{i}\tilde{\theta}^{k}, \qquad i=1,2, \cdots, n,
		\end{equation}
		where  $\tilde{\theta}=1-1/\kappa$ and
		\begin{equation*}
			\left\{\begin{array}{l}
				\widetilde{C}_{1}=\lvert g_{0}^{(1)}\rvert; \\
				\widetilde{C}_{i}=\displaystyle
				\max\left\{
				\lvert g_{0}^{(i)}\rvert,\frac{\lvert g_{1}^{(i)}\rvert}{\tilde{\theta}},
				\cdots,\frac{\lvert g_{r}^{(i)}\rvert}{\tilde{\theta}^{r}}, \frac{\tilde{\sigma}_{i}^{r+1}}{\tilde{\theta}^{r+1}\psi(\lambda_{i})}\sqrt{\sum_{j=1}^{i-1} \psi^{2}(\lambda_{j})\widetilde{C}_{j}^{2}}
				\right\},\\
				\quad i=2,3, \cdots, n,
			\end{array}\right.
		\end{equation*} 
		with  $\tilde{\sigma}_{i}=\max \left\{\frac{\lambda_{i}}{\lambda_{1}}-1,1-\frac{\lambda_{i}}{\lambda_{n}}\right\}$.
		
		\item[(ii)] If Property B is satisfied by replacing $v(k)$ in \eqref{alphack} with $v_i(k)\in\{k,k-1,\cdots,\max\{k-m_i+1,0\}\}$ for all $i=1,\cdots,s$ where $m_i$ and $s$ are positive integers, then either $g_{k} = 0$ for some finite $k$ or the sequence  $\left\{\|g_{k}\|\right\}$ converges to zero $R$-linearly in the sense of \eqref{Citheta3} with $m=\max\{m_1,\cdots,m_s\}$ and 
		\begin{equation*}
			\left\{\begin{array}{l}
				\widetilde{C}_{1}=\lvert g_{0}^{(1)}\rvert; \\
				\widetilde{C}_{i}=\displaystyle
				\max\left\{
				\lvert g_{0}^{(i)}\rvert,\frac{\lvert g_{1}^{(i)}\rvert}{\tilde{\theta}}, \cdots,\frac{\lvert g_{m-1}^{(i)}\rvert}{\tilde{\theta}^{m-1}},\frac{\max \{\tilde{\sigma}_{i}, \tilde{\sigma}_{i}^{m}\}}{\tilde{\theta}^{m}\psi(\lambda_{i})}\sqrt{\sum_{j=1}^{i-1} \psi^{2}(\lambda_{j})\widetilde{C}_{j}^{2}}
				\right\},\\
				\quad i=2,3, \cdots, n.
			\end{array}\right.
		\end{equation*}
	\end{enumerate}
\end{theorem}

\begin{remark}
	Notice that $\alpha_{k}^{BB1}$ satisfies Property B with $v(k)=k-1$, $\psi(A)=I$ and $M_{1}=\lambda_{n}$. So, the result in Theorem 1 of \cite{li2021on} is a special case of (i) in Theorem \ref{converth2}. 	
	The stepsize $\alpha_{k}^{GMR}$ in \eqref{gmr} satisfies Property B with $\psi(A)=A^{\rho(k)/2}$ for some  real number $\rho(k)\geqslant 0$. Hence, the results in Lemma 2.4 of \cite{huang2021On} will be recovered from (ii) of Theorem \ref{converth2} by setting $s=1$ and $m_1=m$ but the last term in $\widetilde{C}_{i}$ is slightly different. It is worth to mention that since $\lambda_{j}^{\rho(k)/2}\leqslant \lambda_{i}^{\rho(k)/2}$ when $\lambda_{j}\leqslant \lambda_{i}$ the quantities $\widetilde{C}_{i}$, $i=2,3, \cdots, n$, here are smaller than or equal to those in \cite{huang2021On}.
\end{remark}




\begin{remark}
	From the above analysis, it seems that for two gradient methods corresponding to the same $v(k)$ and $M_{1}$ the one with smaller $\psi(\lambda_{j})/\psi(\lambda_{i})$  in $\widetilde{C}_{i}$ will have faster convergence rate. However, this is not necessarily true because the components of $g_k$ obtained by the two methods may be quite different. For example, both the BB1 and BB2 methods correspond to $v(k)=k-1$ and $M_{1}=\lambda_{n}$ but to $\psi(A)=I$ and $\psi(A)=A$, respectively, which implies a smaller $\psi(\lambda_{j})/\psi(\lambda_{i})$ for the BB2 method.
	Nevertheless, there is no theoretical evidence  in support of a faster convergence rate for the BB2 method and practical experience often prefers the BB1 method, see \cite{dai2005asymptotic} for example.
\end{remark}


\section{Some discussions}\label{sec:3}
We introduced Property B for gradient methods and proved that any gradient method with stepsizes satisfying this property has $R$-linear convergence rate $1-\lambda_{1}/M_1$ where $M_1$ is the upper bound of the inverse stepsize. For a large collection of gradient methods, we refined the rate to $1-1/\kappa$. An important advantage of Property B is that it can easily be checked. Our results extend the one in \cite{li2021on} and partly improve those in \cite{huang2021On}.


%


With the notation $G(k, l)=\sum_{i=1}^{l}\left(\psi(\lambda_{i})g_{k}^{(i)}\right)^{2}$, we can generalize Property A in \cite{dai2003alternate} as the following Property GA.

\noindent\textbf{Property~GA:} If there exist an integer $m$ and  positive constants $M_{1} \geqslant  \lambda_{1} $ and $M_{2}$ such that\\
(i) $\lambda_{1} \leqslant  \alpha_{k}^{-1} \leqslant  M_{1}$;\\
(ii) for any integer $ l \in [1,n-1] $ and real number $\varepsilon >0 $, if $G(k-j,l) \leqslant  \varepsilon$ and $\left(\psi(\lambda_{l+1})g_{k-j}^{(l+1)}\right)^{2} \geqslant  M_{2}\varepsilon$ hold for $j \in [0,\min\{k,m\}-1]$, then $\alpha_{k}^{-1} \geqslant  \frac{2}{3}\lambda_{l+1}$.	


Clearly, Property GA reduces to Property A when $\psi(A)=I$. Moreover, it is not difficult to check that  Property B is a special case of Property GA. 
By similar arguments as the proof of Theorem \ref{converth}, we can establish $R$-linear convergence rate as \eqref{Citheta} for any gradient method using stepsizes satisfying Property GA but with
\begin{equation*}
	\left\{\begin{array}{l}
		C_{1}=\lvert g_{0}^{(1)}\rvert; \\
		C_{i}=\displaystyle
		\max\left\{
		\lvert g_{0}^{(i)}\rvert,\frac{\lvert g_{1}^{(i)}\rvert}{\theta}, \cdots,\frac{\lvert g _{m-1}^{(i)}\rvert}{\theta^{m-1}},\frac{\max \{\sigma_{i}, \sigma_{i}^{m}\}}{\theta^{m}\psi(\lambda_{i})}\sqrt{\sum_{j=1}^{i-1}M_{2}\psi^{2}(\lambda_{j}) C_{j}^{2}}
		\right\},\\
		\quad i=2,3, \cdots, n,
	\end{array}\right.
\end{equation*}
and $\theta=\max\{\frac{1}{2},1-\frac{\lambda_{1}}{M_1}\}$ which is bounded by $\max\{\frac{1}{2},1-\frac{\lambda_{1}}{M}\}$ with $M=\lim\sup_{k\rightarrow\infty}\alpha_k^{-1}$. This extends the results in Lemma 2.2 and Theorem 2.1 of \cite{huang2021On}. Other results in the former section could be obtained for Property GA as well.

All the above results on gradients can be applied to get $R$-linear convergence for iterations and function values. In fact, denote by $x^*=A^{-1}b$ the solution of \eqref{1a} and $e_k=x_k-x^*$. Then we have
\begin{equation*}
	e_k=x_k-x^*=A^{-1}(Ax_k-b)=A^{-1}g_k
\end{equation*}
and 
\begin{equation*}
	f(x_k)=f(x^*)+\frac{1}{2}e_k^TAe_k.
\end{equation*}
By applying \eqref{Citheta} we obtain
\begin{equation*}
	\lvert e_k^{(i)}\rvert=\lvert\lambda_i^{-1}g_k^{(i)}\rvert\leqslant  \lambda_i^{-1}C_{i}\theta^{k}, \qquad i=1,2, \cdots, n,
\end{equation*}
which yield
\begin{align*}
	f(x_k)-f(x^*)&=\frac{1}{2}e_k^TAe_k
	=\frac{1}{2}\sum_{i=1}^{n}\lambda_i(e_k^{(i)})^2\\
	&\leqslant \frac{1}{2}\sum_{i=1}^{n}\lambda_i^{-1}C_{i}^2\theta^{2k}
	=\frac{1}{2}\|C\|^2\theta^{2k},
\end{align*}
where $C=(\lambda_1^{-1/2}C_{1},\cdots,\lambda_n^{-1/2}C_{n})^T$. In the same way, we are able to establish $R$-linear convergence for iterations and function values from \eqref{Citheta3}.

Notice that the upper bound on the right hand side of \eqref{ineqlin1} can be refined as $\hat{\theta}=\frac{\kappa-1}{\kappa+1}$. From the proof of Theorem \ref{converth} we know that the convergence rate in Theorem \ref{converth2} will be improved to $\hat{\theta}$ if $\lvert g_{k}^{(1)}\rvert\leqslant  \hat{\theta}\lvert g_{k-1}^{(1)}\rvert$ holds. Nevertheless, the above inequality requires  $1-\alpha_{k}\lambda_1\leqslant  \hat{\theta}$, i.e. $\alpha_{k}\geqslant \frac{2}{\lambda_1+\lambda_n}$ which always yields the rate $\hat{\theta}$. The authors of \cite{huang2021On,li2021on} have constructed examples for the method \eqref{gmr} and the BB1 method such that $\alpha_{k}\equiv\frac{2}{\lambda_1+\lambda_n}$ to achieve the rate $\hat{\theta}$, respectively. Recall that both BB1 and \eqref{gmr} are special cases of Property B. Actually, for a given method satisfies Property B, we can construct examples as \cite{huang2021On,li2021on} such that $\alpha_{k}\equiv\frac{2}{\lambda_1+\lambda_n}$. It is worth to mention that the convergence results obtained here and in \cite{huang2021On,li2021on} are componentwise. One interesting question is whether we could get better convergence rate for gradient norm or objective values under reasonable conditions. This will be our future work.

\begin{acknowledgements}
This work was supported by the National Natural Science Foundation of China  (Grant No. 11701137) and Natural Science Foundation of Hebei Province (Grant No. A2021202010).
\end{acknowledgements}


\begin{thebibliography}{10}
	\providecommand{\url}[1]{#1}
	\csname url@samestyle\endcsname
	\providecommand{\newblock}{\relax}
	\providecommand{\bibinfo}[2]{#2}
	\providecommand{\BIBentrySTDinterwordspacing}{\spaceskip=0pt\relax}
	\providecommand{\BIBentryALTinterwordstretchfactor}{4}
	\providecommand{\BIBentryALTinterwordspacing}{\spaceskip=\fontdimen2\font plus
		\BIBentryALTinterwordstretchfactor\fontdimen3\font minus
		\fontdimen4\font\relax}
	\providecommand{\BIBforeignlanguage}[2]{{%
			\expandafter\ifx\csname l@#1\endcsname\relax
			\typeout{** WARNING: IEEEtran.bst: No hyphenation pattern has been}%
			\typeout{** loaded for the language `#1'. Using the pattern for}%
			\typeout{** the default language instead.}%
			\else
			\language=\csname l@#1\endcsname
			\fi
			#2}}
	\providecommand{\BIBdecl}{\relax}
	\BIBdecl
	
	
	\bibitem{cauchy1847methode}
	Cauchy, A.: M{\'e}thode g{\'e}n{\'e}rale pour la r{\'e}solution des systemes
	di'{\'e}quations simultan{\'e}es.
	\newblock Comp. Rend. Sci. Paris \textbf{25}, 536--538 (1847)
	
	\bibitem{akaike1959successive}
	Akaike, H.: On a successive transformation of probability distribution and its
	application to the analysis of the optimum gradient method.
	\newblock Ann. Inst. Stat. Math. \textbf{11}(1), 1--16 (1959)
	
	\bibitem{forsythe1968asymptotic}
	Forsythe, G.E.: On the asymptotic directions of the $s$-dimensional optimum
	gradient method.
	\newblock Numer. Math. \textbf{11}(1), 57--76 (1968)
	
	
	
	
	\bibitem{dai2003altermin}
	Dai, Y.H., Yuan, Y.X.: Alternate minimization gradient method.
	\newblock IMA J. Numer. Anal. \textbf{23}(3), 377--393 (2003)
	
	\bibitem{huang2019asymptotic}
	Huang, Y.K., Dai, Y.H., Liu, X.W. et~al.: On the asymptotic convergence and acceleration of gradient methods. J. Sci. Comput. \textbf{90}, 7, (2022).
	
	\bibitem{dai2006new2}
	Dai, Y.H., Yang, X.: A new gradient method with an optimal stepsize property.
	\newblock Comp. Optim. Appl. \textbf{33}(1), 73--88 (2006)
	
	\bibitem{elman1994inexact}
	Elman, H.C., Golub, G.H.: Inexact and preconditioned Uzawa algorithms for
	saddle point problems.
	\newblock SIAM J. Numer. Anal. \textbf{31}(6), 1645--1661 (1994)
	
	\bibitem{dai2005analysis}
	Dai, Y.H., Yuan, Y.X.: Analysis of monotone gradient methods.
	\newblock J. Ind. Mang. Optim. \textbf{1}(2), 181 (2005)
	
	\bibitem{yuan2006new}
	Yuan, Y.X.: A new stepsize for the steepest descent method.
	\newblock J. Comput. Math. \textbf{24}(2), 149--156 (2006)
	
	\bibitem{yuan2008step}
	Yuan, Y.X.: Step-sizes for the gradient method.
	\newblock AMS/IP Stud. Adv. Math.  \textbf{42}(2), 785--796
	(2008)
	
	\bibitem{barzilai1988two}
	Barzilai, J., Borwein, J.M.: Two-point step size gradient methods.
	\newblock IMA J. Numer. Anal. \textbf{8}(1), 141--148 (1988)
	
	\bibitem{dai2013new}
	Dai, Y.H.: A new analysis on the Barzilai--Borwein gradient method.
	\newblock J. Oper. Res. Soc. China \textbf{2}(1), 187--198 (2013)
	
	
	\bibitem{dai2005asymptotic}
	Dai, Y.H., Fletcher, R.: On the asymptotic behaviour of some new gradient
	methods.
	\newblock Math. Program. \textbf{103}(3), 541--559 (2005)
	
	\bibitem{raydan1993barzilai}
	Raydan, M.: On the Barzilai--Borwein choice of steplength for the gradient
	method.
	\newblock IMA J. Numer. Anal. \textbf{13}(3), 321--326 (1993)
	
	\bibitem{dai2002r}
	Dai, Y.H., Liao, L.Z.: $R$-linear convergence of the Barzilai--Borwein
	gradient method.
	\newblock IMA J. Numer. Anal. \textbf{22}(1), 1--10 (2002)
	
	
	\bibitem{li2021on}
	Li, D.W., Sun, R.Y.: On a faster $R$-Linear convergence rate of the Barzilai--Borwein method, \newblock arXiv preprint arXiv:2101.00205v2, (2021)
	
	\bibitem{dai2015positive}
	Dai, Y.H., Al-Baali, M., Yang, X.: A positive Barzilai--Borwein-like stepsize
	and an extension for symmetric linear systems.
	\newblock In: Numerical Analysis and Optimization, pp. 59--75. Springer (2015)
	
	
	\bibitem{BDH2019stabilized}
	Burdakov, O., Dai, Y. H., Huang, N.: Stabilized Barzilai-Borwein Method, J. Comp. Math., \textbf{37}(6), 916--936 (2019)
	
	
	\bibitem{dai2019family}
	Dai, Y.H., Huang, Y.K., Liu, X.W.: A family of spectral gradient methods for
	optimization.
	\newblock Comp.
	Optim. Appl. \textbf{74}(1), 43--65
	(2019)
	
	\bibitem{di2018steplength}
	Di~Serafino, D., Ruggiero, V., Toraldo, G., Zanni, L.: On the steplength
	selection in gradient methods for unconstrained optimization.
	\newblock Appl. Math. Comput. \textbf{318}, 176--195 (2018)
	
	
	\bibitem{fletcher2005barzilai}
	Fletcher, R.: On the Barzilai--Borwein method.
	\newblock  In: Optimization and Control with Applications, pp. 235--256. Springer, New York (2005)
	
	\bibitem{huang2020equipping}
	Huang, Y.K., Dai, Y.H., Liu, X.W.: Equipping the Barzilai--Borwein method with the two
	dimensional quadratic termination property.
	\newblock  SIAM J. Optim. \textbf{31}(4), 3068--3096 (2021)
	
	\bibitem{raydan1997barzilai}
	Raydan, M.: The Barzilai--Borwein gradient method for the large scale
	unconstrained minimization problem.
	\newblock SIAM J. Optim. \textbf{7}(1), 26--33 (1997)
	
	\bibitem{birgin2000nonmonotone}
	Birgin, E.G., Mart{\'\i}nez, J.M., Raydan, M.: Nonmonotone spectral projected
	gradient methods on convex sets.
	\newblock SIAM J. Optim. \textbf{10}(4), 1196--1211 (2000)		
	
	\bibitem{crisci2020spectral}
	Crisci, S., Porta, F., Ruggiero, V., Zanni, L.: Spectral properties of
	Barzilai--Borwein rules in solving singly linearly constrained optimization
	problems subject to lower and upper bounds.
	\newblock  SIAM J. Optim.  \textbf{30}(2), 1300--1326 (2020)
	
	
	\bibitem{dai2005projected}
	Dai, Y.H., Fletcher, R.: Projected Barzilai--Borwein methods for large-scale
	box-constrained quadratic programming.
	\newblock Numer. Math. \textbf{100}(1), 21--47 (2005)
	
	
	\bibitem{huang2016smoothing}
	Huang, Y.K., Liu, H.: Smoothing projected Barzilai--Borwein method for
	constrained non-lipschitz optimization.
	\newblock Comp. Optim. Appl. \textbf{65}(3), 671--698 (2016)
	
	\bibitem{huang2020gradient}
	Huang, Y.K., Dai, Y.H., Liu, X.W., Zhang, H.: Gradient methods exploiting
	spectral properties.
	\newblock Optimi. Method Softw. \textbf{35}(4), 681--705 (2020)
	
	\bibitem{huang2015quadratic}
	Huang, Y.K., Liu, H., Zhou, S.: Quadratic regularization projected
	Barzilai--Borwein method for nonnegative matrix factorization.
	\newblock Data Min. Knowl. Disc. \textbf{29}(6), 1665--1684 (2015)
	
	\bibitem{jiang2013feasible}
	Jiang, B., Dai, Y.H.: Feasible Barzilai--Borwein-like methods for extreme
	symmetric eigenvalue problems.
	\newblock Optim. Method Softw. \textbf{28}(4), 756--784 (2013)
	
	\bibitem{tan2016barzilai}
	Tan, C., Ma, S., Dai, Y.H., Qian, Y.: Barzilai--Borwein step size for
	stochastic gradient descent.
	\newblock In: Advances in Neural Information Processing Systems, pp. 685--693
	(2016)
	
	\bibitem{wang2007projected}
	Wang, Y., Ma, S.: Projected Barzilai--Borwein method for large-scale
	nonnegative image restoration.
	\newblock Inverse Probl. Sci. En. \textbf{15}(6), 559--583 (2007)
	
	\bibitem{de2014efficient}
	De~Asmundis, R., Di~Serafino, D., Hager, W.W., Toraldo, G., Zhang, H.: An
	efficient gradient method using the Yuan steplength.
	\newblock Comp. Optim. Appl. \textbf{59}(3), 541--563 (2014)
	
	\bibitem{frassoldati2008new}
	Frassoldati, G., Zanni, L., Zanghirati, G.: New adaptive stepsize selections in
	gradient methods.
	\newblock J. Ind. Mang. Optim. \textbf{4}(2), 299--312 (2008)
	
	\bibitem{sun2020new}
	Sun, C., Liu, J.P.: New stepsizes for the gradient method.
	\newblock Optim. Lett. \textbf{14}(7), 1943--1955 (2020)
	
	
	\bibitem{zhou2006gradient}
	Zhou, B., Gao, L., Dai, Y.H.: Gradient methods with adaptive step-sizes.
	\newblock Comp. Optim. Appl. \textbf{35}(1), 69--86 (2006)
	
	\bibitem{dai2003alternate}
	Dai, Y.H.: Alternate step gradient method.
	\newblock Optimization \textbf{52}(4-5), 395--415 (2003)
	
	
	\bibitem{huang2021On}
	Huang, N.: On $R$-linear convergence analysis for a class of gradient methods.
	\newblock Comput. Optim. Appl. \textbf{81}(1), 161--177 (2022).
	
	\bibitem{friedlander1998gradient}
	Friedlander, A., Mart{\'\i}nez, J.M., Molina, B., and Raydan, M.:
	Gradient method with retards and generalizations.
	\newblock SIAM J. Numer. Anal. \textbf{36}(1), 275--289 (1998)
	
	
	\bibitem{huang2020acceleration}
	Huang, Y.K., Dai, Y.H., Liu, X.W., Zhang, H.: On the acceleration of the
	Barzilai--Borwein method.
	\newblock Comp. Optim. Appl. \textbf{81}(3), 717--740 (2022)
	
	
	
	
	\bibitem{zou2021fast}
	Zou, Q., Magoul{\`e}s, F.: Fast gradient methods with alignment for symmetric
	linear systems without using Cauchy step.
	\newblock J. Comput. Math.  \textbf{381}, 113033
	(2021)
	
	
	\bibitem{dai2006cyclic}
	Dai, Y.H., Hager, W.W., Schittkowski, K., Zhang, H.: The cyclic
	Barzilai--Borwein method for unconstrained optimization.
	\newblock IMA J. Numer. Anal. \textbf{26}(3), 604--627 (2006)
	
	
	
\end{thebibliography}

\end{document}